\documentclass[10pt, reqno, final]{amsart}

\usepackage{caption}
\usepackage{float}
\usepackage[T1]{fontenc}
\usepackage[utf8]{inputenc}
\usepackage[english]{babel}
\usepackage{graphicx}
\usepackage{float}
\usepackage{amsaddr}
\usepackage{fancyhdr}

\usepackage{amsmath}
\usepackage{amssymb}
\usepackage{amsthm}

\usepackage[backend=biber,url=false,doi=false,isbn=false,firstinits=true, hyperref=true,maxbibnames=99]{biblatex}
\theoremstyle{plain}

\newtheorem{thm}{Theorem}[section]

\newtheorem{cor}[thm]{Corollary}

\newtheorem{lem}[thm]{Lemma}		

\newtheorem{prop}[thm]{Proposition}

\theoremstyle{definition}
\newtheorem{defn}[thm]{Definition}

\theoremstyle{remark}
\newtheorem{rem}[thm]{Remark}

\newtheorem*{assumption*}{\assumptionnumber}
\providecommand{\assumptionnumber}{}
\makeatletter
\newenvironment{assumption}[2]
{%
	\renewcommand{\assumptionnumber}{\textup{(H#1)}}%
	\begin{assumption*}%
		\protected@edef\@currentlabel{#1}%
	}
	{%
	\end{assumption*}
}

\numberwithin{equation}{section}
\usepackage{mathtools}
\usepackage{mathrsfs}
\usepackage{eucal}
\usepackage{braket}
\usepackage{a4wide}
\usepackage{booktabs}		
\usepackage{multirow}		
\usepackage{caption}		
\usepackage{rotating}
\usepackage{subfig}
\usepackage{varioref}		
\usepackage{footmisc}

\usepackage{epsfig}
\usepackage{subfig}
\usepackage{enumerate}
\usepackage[colorlinks=true, linkcolor=blue, citecolor=blue, 
urlcolor=blue]{hyperref}		
\mathtoolsset{showonlyrefs}
\input xy
\xyoption{all}
\usepackage{geometry}
\newcommand{\R}{\mathbb{R}}
\newcommand{\Com}{\mathbb{C}}
\newcommand{\Q}{\mathbb{Q}}
\newcommand{\Z}{\mathbb{Z}}	
\newcommand{\N}{\mathbb{N}}	
\renewcommand{\P}{\mathscr{P}}
\newcommand{\Lagr}{\mathrm{Lag}}
\newcommand{\Bsym}{\mathrm{B}_{\textup{sym}}}
\newcommand{\norm}[1]{\left\| #1 \right\|}		
\newcommand{\cfsa}{\mathscr{CF}^{sa}}
\newcommand{\im}{\mathrm{rge}\,}
\newcommand{\Mat}{\mathrm{Mat}\,}
\newcommand{\Sp}{\mathrm{Sp}\,}
\newcommand{\Sym}{\mathrm{Sym}\,}
\newcommand{\Lin}{\mathscr{L}}
\newcommand{\GL}{\mathrm{GL}}
\newcommand{\Id}{I}
\newcommand{\br}{~\\}
\newcommand{\re}{\mathrm{\mathfrak{R}}}
\newcommand{\imaginary }{\mathrm{\mathfrak{I}}}
\newcommand{\ev}{\mathrm{ev}}
\newcommand{\Sf}[1]{\mathrm{Sf\,}(#1)}
\newcommand{\Inp}[1]{\left\langle#1\right\rangle}
\DeclareMathOperator{\diag}{diag}
\newcommand{\trasp}[1]{{#1}^\mathsf{T}}	
\DeclareMathOperator{\sgn}{sgn}
\DeclareMathOperator{\ind}{ind}
\DeclareMathOperator{\iCLM}{\iota^{\scriptscriptstyle{\mathrm{CLM}}}}
\DeclareMathOperator{\iMor}{m^{--}}
\DeclareMathOperator{\coiMor}{\iota^{+}}
\DeclareMathOperator{\irel}{i_{\textup{rel}}}
\DeclareMathOperator{\iRS}{\mu^{\textup{RS}}}
\DeclareMathOperator{\ispec}{\iota_{\textup{spec}}}
\DeclareMathOperator{\igeo}{\iota}
\DeclareMathOperator{\igeob}{\iota_b}
\DeclareMathOperator{\dif}{d}
\DeclareMathOperator{\dom}{dom}
\DeclareMathOperator{\codim}{codim}
\renewcommand{\leq}{\leqslant}
\renewcommand{\geq}{\geqslant}
\renewcommand{\tilde}{\widetilde}
\renewcommand{\=}{\coloneqq}
\newcommand{\eq}{\eqqcolon}
\newcommand{\todo}[1]{\fbox{\large ** TO DO: #1 **}}
\newcommand{\n}[1]{{\bf #1}}

\title{Index theory for traveling waves in reaction diffusion systems with skew gradient structure }

\author{ Qin Xing} 
\thanks{Partially supported by NSFC  No.11425105}
\address{School of Mathematics,
	Shandong University,
	Jinan, Shandong, 
	The People's Republic of China }
\email{xingqin@mail.sdu.edu.cn}
\date{}
\begin{document}
\maketitle	
\thispagestyle{fancy}
\lhead{} 
\chead{} 
\rhead{} 
\lfoot{} 
\cfoot{\thepage} 
\pagestyle{fancy}
\rfoot{}
\renewcommand{\headrulewidth}{0pt} 
\renewcommand{\footrulewidth}{0pt} 
\begin{abstract}
	
	A unified geometric approach for the stability analysis of traveling pulse solutions for reaction diffusion equations with skew-gradient structure has been established in a previous paper \cite{cornwell2019openinga}, but essentially no results have been found in the case of traveling front solutions. In this work, we will bridge this gap. For such cases, a Maslov index of the traveling wave is well-defined, and we will show how it can be used to provide the spectral information of the waves. As an application, we use the same index providing the exact number of unstable eigenvalues of the traveling front solutions of FitzHugh-Nagumo equations.

	\vskip0.2truecm
	\noindent
	\vskip0.1truecm
	\noindent
	\textbf{Keywords:} heteroclinic orbits; stability of traveling wave; Maslov index; spectral flow; FitzHugh–Nagumo equations
	\vskip0.5truecm 
\end{abstract}	

	\section{Introduction and main results}


The setting of this paper is the systems of reaction-diffusion equations of the form	
	\begin{equation}
	w_t=w_{xx}+Q\nabla F(w),\label{eq:r.e.q.}
	\end{equation}
	where $x,t\in\R$ are space and time, respectively, $w\in\R^{n}$ and $\nabla F$ are the gradients of the function $F:\R^{n}\rightarrow\R$,	and $Q\in \Mat(n,\R)$ has the form $Q=\begin{bmatrix}
	\Id_r&0\\0&-\Id_{n-r}
	\end{bmatrix}$, where $\Mat(n,\R)$ is the set of all $n\times n$ matrices. Such a system is of the activator-inhibitor type and will be referred to as a skew-gradient system \cite{yanagida2002minimaximizers}. 
	 A traveling wave solution $w^*$ of one variable $\xi=x-ct$ is a solution to   \eqref{eq:r.e.q.} with asymptotic behavior \[w_{\pm} \text{ as } \xi\rightarrow\pm\infty.\] 
	Here, $w_{\pm}$ are the constant equilibria of \eqref{eq:r.e.q.}: that is, $\nabla F(w_{\pm})=0$. If $w_+=w_- $, then $w^*$ is called a pulse (and a front otherwise). In this paper, we are focused primarily on fronts. We will always consider the case of $c>0$: otherwise, the direction of motion can be reversed. 
	
		Traveling waves being solutions of \eqref{eq:r.e.q.} is an important subjection in dynamical systems, and many aspects of traveling waves have been studied in extensive works\cite{cornwell2018existence,chen2018front}. One of the key issues of concern is whether a steady state is stable with respect to disturbances under initial conditions since this directly determines whether it can be observed in nature. Motivated by \cite{cornwell2019openinga}, in this paper, we establish a unified geometric approach for the stability analysis of traveling front solutions    for \eqref{eq:r.e.q.}.
	Written in a moving frame, a traveling front solution $w^*$ of (\ref{eq:r.e.q.}) can be regarded as a heteroclinic solution $w^*$ of the following equation:
	\begin{equation}
	\begin{cases}
	w_{\xi\xi}+cw_\xi+Q\nabla F(w) & =0,\\
	\lim_{\xi\rightarrow\pm\infty}w(\xi)=w_{\pm}.
	\end{cases}\label{eq:ODE from R.E.Q.}
	\end{equation}
	The stability analysis is directly related to the spectral information of the operator $ L:=\frac{\dif^2}{\dif \xi^2}+c\frac{\dif}{\dif \xi}+QB(\xi)$ by linearizing \eqref{eq:ODE from R.E.Q.} along $w^*$, where $B(\xi)=\nabla^2F(w^*)$; note that the matrices $B_\pm\=\lim\limits_{\xi\to\pm\infty}B(\xi)$ are well-defined. Moreover, the previous description guarantees that there exists $C>$, such that \begin{align}\label{eq:constant C}
	\Inp{QB(\xi)v,v}\leq C|v|^2 \text{ for all } (\xi,v)\in\R\times \R^n.
	\end{align}
	
	 For $z\in\Com $, denote by $\Re z$ and $\imaginary  z$ the real and imaginary parts of $z$, respectively. $\R^+\=(0,+\infty)$,  $\R^-\=(-\infty,0)$, $ \Com^+\=\{z\in\Com|\Re z> 0 \}$ and $  \Com^-\=\{z\in\Com|\Re z< 0 \}$. The pair $(\R^n,\langle,\rangle)$ denotes the $n$-dimensional Euclidean space. Moreover, $\overline{\#}$	denotes the closure of the set $\#$. 
	 
	 Since a wave solution of  \eqref{eq:r.e.q.} possesses a translation invariance property, it is said to be nondegenerate if zero is a simple eigenvalue of $L$.
	 
	\begin{defn}\label{def:stable}
		A nondegenerate wave solution of  \eqref{eq:r.e.q.} is spectrally stable if all the nonzero eigenvalues of $L$ are in $\Com^-$.
	\end{defn}


	In this paper, we study the following special case.
	
\begin{assumption}{1}{2}\label{as:1}
	Suppose that $\sigma\left(Q B_{\pm}\right)\subset \Com^{-} $.
\end{assumption}
	Thus, $w_+$ and $w_-$ are both stable equilibria of \eqref{eq:r.e.q.}.



	From now on, we focus on studying the eigenvalue problem \begin{equation}\label{eq:eigen.p.}
	L\phi=\lambda\phi.
	\end{equation} Denoted by $\sigma_p(L)$ is the set of isolated eigenvalues with finite multiplicity, and $\sigma_{ess}(L)=\sigma(L)\backslash\sigma_p(L)$ is the essential spectrum of $L$. It is known (Cf. Lemma \ref{lem:eigenvalue order of A}) that   $\sigma_{ess} (L)\subset \Com^-$. As was commonly performed in \cite{cornwell2019openinga}, we set $y=\begin{bmatrix}
	\dot{\phi}\\\phi
	\end{bmatrix}$ to convert  \eqref{eq:eigen.p.} to \begin{equation}\label{eq:ODE}
	\dot{y}=A_\lambda (\xi)y,
	\end{equation}
	where $A_\lambda(\xi)=\begin{bmatrix}
	-c&\lambda \Id-Q B(\xi)\\\Id &0
	\end{bmatrix}$, and there are well-defined matrices $A_\lambda(\pm\infty)=\lim\limits_{\xi\to\pm\infty}A_\lambda(\xi)$. Throughout this paper, the dots denote differentiation with respect to $\xi$.

	For any $M\in \Mat(\R,\R^n)$, denote by $V^+(M)$ ($V^-(M)$) the positive (negative) spectral space corresponding to the eigenvalues of $M$ having positive (negative) real parts.
	 We provide the next two conditions:
	 \begin{assumption}{2}{2}\label{as:3}
$(Q B_{\pm}v,v)<0$, for all $v\in V^-(Q)\backslash\{0\}$,
	 \end{assumption}
 \begin{assumption}{2'}{2}\label{as:4}
$(Q B_{\pm}v,v)<0$, for all $v\in\R^n\backslash\{0\}$.
 \end{assumption}
	 
Letting $B\in \Mat(n,\R)$, in the following remark, we will explain the relationship between (H\ref{as:1}) and (H\ref{as:4}).	
\begin{rem}\label{rem:condition}
(1)	We claim that (H\ref{as:4}) implies that (H\ref{as:1}). Let $\lambda$ be the eigenvalue of $QB$ with eigenvector $v$. If $\lambda\in \R$, then we have that $\lambda|v|^2=\lambda\Inp{v,v}=\Inp{QBv,v}<0$, which implies that $\lambda<0$. If $\lambda\in \Com\backslash\R$, then $\bar \lambda$ is also the eigenvalue of $QB$ with eigenvector $\bar v$, where $\bar \lambda$ and $\bar v$ denote the complex conjugate of $\lambda$ and $v$, respectively. Now, let $\lambda=a+bi$ and $v=u+iw$: then, we have that \begin{align}
\lambda(|u|^2+|w|^2)=\lambda\Inp{v,v}=\Inp{QBv,v}=\Inp{QBu,u}-i\Inp{QBu,w}+i\Inp{QBw,u}+\Inp{QBw,w}.
\end{align} Similarly, we have that \[\bar \lambda(|u|^2+|w|^2)=\Inp{QBu,u}+i\Inp{QBu,w}-i\Inp{QBw,u}+\Inp{QBw,w}.\] We thus determine that $ a=\frac{1}{2}(\lambda+\bar\lambda)=\frac{\Inp{QBu,u}+\Inp{QBw,w}}{|u|^2+|w|^2}<0.$ This implies that $\lambda$ has a negative real part.

(2) We remark that (H\ref{as:4}) cannot be derived from (H\ref{as:1}). We consider $B\in \Mat(2,\R)$ with the form $B=\begin{bmatrix}
1&\frac{\sqrt{35}}{4}\\\frac{\sqrt{35}}{4}&2
\end{bmatrix}$ and $Q=\begin{bmatrix}
1&0\\0&-1
\end{bmatrix}$. It is easy to determine that $\sigma(QB)\subset\Com^-$, but $\Inp{QB\begin{bmatrix}
	1\\0
	\end{bmatrix},\begin{bmatrix}
	1\\0
	\end{bmatrix}}=1>0$.

(3) Supposing that $ \sigma(QB)\subset \Com^-$, we claim that there exists a nonsingular matrix $T\in \Mat(n,\R)$, such that \[ \langle  T^{-1}QBTv,v\rangle<0,\ \forall v\in \R^n\backslash\{0\} .\] We only consider a special case, where the general case is analogous. Letting $QB\in\Mat(6,\R)$, suppose that $\lambda_1\in\R$ and $\lambda_2=a+ib$ are generalized eigenvalues of $QB$ and have the same algebraic multiplicity 2: then, for $\epsilon <-\max\{\lambda_1,a\}$, there exists $T_\epsilon$, such that \[T_\epsilon^{-1}QBT_\epsilon=\begin{bmatrix}
\lambda_1&\epsilon&0&0&0&0\\0&\lambda_1&0&0&0&0\\0&0&a&-b&\epsilon&0\\0&0&b&a&0&\epsilon\\0&0&0&0&a&-b\\0&0&0&0&b&a
\end{bmatrix}.\] By a simple calculation, we have that \begin{align}
\langle T_\epsilon^{-1}QBT_\epsilon v,v \rangle&=\lambda_1 v_1^2+\lambda_1 v^2_2+\epsilon v_1v_2+av^2_3+av^2_4+av^2_5+av^2_6+\epsilon v_3v_5+\epsilon v_4v_6\\ &\leq \sum_{i=1}^{2}(\lambda_1+\frac{1}{2}\epsilon) v_i^2+\sum_{j=3}^6(a+\frac{1}{2}\epsilon )v^2_j<0.
\end{align}
\end{rem}


 Following Definition \ref{def:stable}, we focus our attention on $\lambda\in \overline{\R^+}$: suppose that $\Phi_{\tau,\lambda}(\xi)$ is the matrix solution of \eqref{eq:ODE}, such that $\Phi_{\tau,\lambda}(\tau)=\Id$.
We recall that the stable and unstable subspaces of \eqref{eq:ODE} are \begin{equation}
E_{\lambda}^{s}(\tau):=\left\{ v\in\R^{2n}|\lim_{\tau\rightarrow+\infty}\Phi_{\tau,\lambda}(\xi)v=0\right\} \ \text{and \ }E_{\lambda}^{u}(\tau):=\left\{ v\in\R^{2n}|\lim_{\tau\rightarrow-\infty}\Phi_{\tau,\lambda}(\xi)v=0\right\}.
\end{equation} Throughout this paper, we abbreviate $E_0^s(\tau)$ and $E^u_0(\tau)$ as $E^s(\tau)$ and $E^u(\tau)$, respectively.

To realize the symplectic structure, we introduce the matrix 
$J\=\begin{bmatrix}
0&-Q\\Q&0
\end{bmatrix}$. 
Since $Q^{2}=\Id$ and $Q^{T}=Q$, it follows that $J^{2}=-\Id$ and $J^{T}=-J$. Therefore, $J$ is a complex structure on $\R^{2n}$, and $\omega(\cdot,\cdot)\=\Inp{J\cdot,\cdot}$ defines a syplectic form on $\R^{2n}$.  $\mathrm{Lag}(n)$ denotes the set of all Lagrangian subspaces of $(\R^{2n},\omega)$.

By invoking \cite[Lemma 3.1]{hu2017index}, Remarks \ref{rm:JB hyperbolic} and \ref{rm:eigenvalue equivalent} imply that the subspaces $E^s_\lambda(\tau)$ and $E^u_\lambda(\tau) $ are both Lagrangian for $(\tau,\lambda)\in\R\times\overline{\R^+}$. Moreover, those results can be obtained from the same discussion in \cite[Theorem 2.3]{cornwell2019openinga}.
 
 We denote by $\P([0,1];\R^{2n})$ the space of all ordered pairs of continuous maps of Lagrangian subspaces $ L:[0,1]\ni t\longmapsto L(t)\=\big(L_1(t),L_2(t)\big)\in\Lagr(n)\times\Lagr(n)$ equipped with the compact-open topology. Following the authors in \cite{cappell1994maslov}, we are in the position to briefly recall the definition of the {\em Maslov index for pairs of Lagrangian subspaces,\/} which will be denoted throughout the paper by the symbol $\iCLM$. Loosely speaking, given the pair $L=(L_1,L_2)\in\P([0,1];\R^{2n})$,
this index counts with signs and multiplicities the number of instants
$t\in[0,1]$ that $L_1(t)\cap L_2(t)\neq\{0\}$. 
\begin{defn}\label{def:maslov index}
	The {\em CLM-index\/} is the unique integer valued function 
	\[
	\iCLM:\P([0,1];\R^{2n})\ni L\longmapsto\iCLM(L;[0,1])\in\Z
	\]
	satisfying Properties I-VI given in \cite[Section 1]{cappell1994maslov}. 
\end{defn}
For the sake of the reader, we list a couple of properties of the $\iCLM$-index that we shall use throughout the paper. 
\begin{itemize}
	\item {\bf (Reversal)\/} Let $L\=(L_1,L_2)\in\P([a,b];\R^{2n})$. Denoting by $\widehat L \in\P([-b, -a];\R^{2n})$ the path traveled in the opposite direction, and by setting $\widehat L\=(L_1(-s),L_2(-s))$,  we obtain 
	\[
	\iCLM(\widehat L;[-b,-a])=-\iCLM(L; [a,b]).
	\]
	\item {\bf (Stratum homotopy relative to the ends)} Given a continuous map\[L:[a,b]\ni s\rightarrow L(s)\in \P([a,b];\R^{2n})\text{ where }L(s)(t)\=(L_1(s,t),L_2(s,t)) \]
	such that $\dim L_1(s,a)\cap L_2(s,a)$ and $\dim L_1(s,b)\cap L_2(s,b)$ are both constant, and then, \[\iCLM(L(0);[a,b])=\iCLM(L(1);[a,b]).\]
\end{itemize}

Moreover, one efficient way to study the Maslov index is via the crossing form introduced by \cite{robbin1993maslov} as follows.

Let $L(t):[0,1]\to \Lagr(n)$ be a smooth curve with $L(0)=L_0$. Let $ W$ be a fixed Lagrangian complement of $L(t)$. For $v\in L_0$ and small $t$, define $w(t)\in W$ by $v+w(t)\in L(t)$. The form $Q(v)=\left.\frac{\dif}{\dif t}\right|_{t=0}\omega(v,w(t))$ is independent of the choice of $W$\cite{robbin1993maslov}. A crossing for $L(t)$ is some $t$ for which $L(t)$ intersects $V\in \Lagr                                                                                       (n)$ nontrivially.  At each crossing, the crossing form is defined as \begin{align}
	\Gamma(L(t),V;t)=\left.Q\right|_{L(t)\cap V}.
\end{align} A crossing is called regular if the crossing form is nondegenerate. For a quadratic form $Q $, we use the notation $\text{sign}(Q)$ for its signature. We also write $m^+(Q)$ and $m^-(Q)$ for the positive and negative indices of inertia of $Q$, so that \[\text{sign}(Q)=m^+(Q)-m^-(Q).\] From \cite{zhu1999maslovtypea}, if the path $L(t)$ has only regular crossing with respect to $V$, then we have that \begin{align}\label{eq:maslov index regular crossing}
	\iCLM(V,L(t),t\in[a,b])=m^+(\Gamma(L(a),V;a))+\sum\limits_{a<t<b}\text{sign}\Gamma(L(t),V;t)-m^-(\Gamma(L(b),V;b)).
\end{align}

We recall the definition of Maslov index for a traveling pulse wave $w^*$ of \eqref{eq:r.e.q.} defined in \cite{cornwell2019openinga}. The only requirement is that $\tau_0$ is large enough so that \begin{align}\label{essential condition}
	E^s(\tau_0)\cap E^u(\tau)=\{0\} \text{ for all } \tau\geq \tau_0.
\end{align}

Moreover, since the transformation used in the conversion of \eqref{eq:eigen.p.} to \eqref{eq:ODE} is different from that in \cite{cornwell2019openinga}, then the crossing form defined in the previous paper \cite{cornwell2019openinga} differs from ours by a negative sign, and so, we provide \cite[Definition 3.6]{cornwell2019openinga} in the following form:
\begin{defn}\label{def:another maslov index}\cite[Definition 3.6]{cornwell2019openinga}Let $w^*$ be a traveling pulse solution of \eqref{eq:r.e.q.}, where the Maslov index of $w^*$ is given by  
	$$\mathbf{Maslov}(w^*)\=\sum\limits_{\tau<\tau_0}\mathrm{sign}\Gamma\left(E^u(\tau),E^s(\tau_0);\tau\right)-m^-(\Gamma\left(E^u(\tau),E^s(\tau_0);\tau_0\right)).$$
\end{defn}

By using the Maslov index, this paper \cite{cornwell2019openinga} provides a unified geometric treatment for the stability analysis of the traveling pulses solution and has been successfully employed in \cite{cornwell2018existence,cornwell2020stability} to obtain some elegant stability results. Motivated by  \cite{cornwell2019openinga}, it is natural to consider the case of the traveling front solutions. We remark that Definition \ref{def:another maslov index} is dependent on the condition \eqref{essential condition}. It is easy to check that this condition may be inefficient for the traveling front solutions, and we sidestep this difficulty by following \cite[Definition 1.2]{hu2017index} and define the Maslov index for the traveling waves as following.

\begin{defn}\label{def:maslov index for traveling wave}
	Let $w^*$ be a traveling wave of \eqref{eq:r.e.q.}, and define the Maslov index of $w^*$ as $$\igeo(w^*)=-\iCLM\left(E^s(\tau),E^u(-\tau);\tau\in\overline{\R^+}\right).$$
\end{defn}

\begin{rem}
	If  Definition \ref{def:another maslov index} is well-defined, then it should be independent of the choice of $\tau_0$. It was shown in \cite{chen2007maslova} that  this definition is independent of $\tau_0$, as long as \eqref{essential condition} is satisfied. However, for the traveling front solution, there may be an intersection point between $E^s(\tau)$ and $E^u(-\tau)$ at $+\infty$. Thus, \eqref{essential condition} may be inefficient. We sidestep this difficulty by following \cite[Definition 1.2]{hu2017index}, and the following proposition shows that  Definition \ref{def:maslov index for traveling wave} generalizes Definition \ref{def:another maslov index} to the front case.
\end{rem}
\begin{prop}\label{pro:change index}
	Letting $w^*$ be a traveling pulse solution of \eqref{eq:r.e.q.} and supposing that (H1) holds, then there exists $\tau_0>0$, such that $E^u(-\infty)\cap E^s(+\infty)=\{0\}$ for all $\tau\geq \tau_0$ and \begin{equation}
		\igeo(w^*)=\mathrm{Maslov }(w^*).
	\end{equation}
\end{prop}
 
Lemmas \ref{lem:nondegenerate condition} and \ref{lem:eigenvalue equivalent} show that the set of nonnegative, real eigenvalues of $L$ is bounded above, 
and then, the spectrum of $L$ in $\overline{\Com^+}$ consists of isolated eigenvalues of finite multiplicity (Cf. page 172 of
\cite{alexander1990topological}), and so, it follows that the quantities \begin{itemize}
		\item $N_+(L)\= \text{the number of real, positive eigenvalues of $L$ counting algebraic multiplicity}$
		\item $\overline N_+(L)\= \text{the number of real, nonnegative eigenvalues of $L$ counting algebraic multiplicity}$
\end{itemize} are both well-defined.

    Now, we introduce a symplectic invariant called triple index (Cf. Definition \ref{eq:the triple index}) which
    can be used to compute the Maslov index, and we denote by $\igeo(L_1,L_2;L_3)$ the triple index for any $L_1,L_2,L_3\in \mathrm{Lag}(n)$.

\begin{thm}\label{thm:main result}
	 If (H\ref{as:1}) and (H\ref{as:3}) hold, then  \begin{equation}
	|\igeo(w^*)+\igeo\left(E^u(-\infty),E^s(+\infty);L_R\right)|\leq \overline N_0(L), \text{ where } L_R=\left\{\left. \begin{bmatrix}
		p\\q
	\end{bmatrix} \right|  p\in V^+(Q)\ \mathrm{and}\ q\in V^-(Q)   \right\}.
	\end{equation}
\end{thm}

\begin{rem}
Under the conditions (H1) and (H2), if \eqref{essential condition} holds, then   Definition \ref{def:another maslov index} is well-defined for the traveling front solution. Using the notation of {\em Maslov box\/} in \cite[Figure 4.1]{cornwell2019openinga}, by the same discussion in the proof of Proposition \ref{pro:change index}, it is easy to prove that the Maslov index contribution along the ``bottom shelf'' $\alpha_4$ is equal to $\igeo\left(E^u(-\infty),E^s(+\infty);L_R\right)$.  Based on this discussion and   Theorem \ref{thm:main result}, the distinction between the pulse and front is nontrivial, and we generalize \cite[Theorem 4.1]{cornwell2019openinga} to the front case.
\end{rem}
Now, we present the central result of this paper.
\begin{thm}\label{thm:central result}
	If (H\ref{as:4}) holds, then we have that \begin{equation}
		|\igeo(w^*)|\leq \overline N_0(L).
	\end{equation}
\end{thm}
\begin{rem}
	 If (H\ref{as:4}) holds, by Lemmas \ref{lem:matrix<0} and \ref{lem:matrix>0}, we have that $E^s(+\infty)\cap E^u(-\infty)=\{0\}$, and this implies that if \eqref{essential condition} holds, then Definition \ref{def:another maslov index} is well-defined. By invoking \eqref{eq:compute maslov index} and Lemma \ref{lem:boundary maslov index = 0}, (H\ref{as:4}) provides a sufficient condition for the Maslov index contribution along the ``bottom shelf'' when $\alpha_4$ is equal to $0$. This can be easily verified for $\lambda=0$, and it  remains true for $\lambda>0$. This is important for practitioners who want to use the Maslov index to prove stability. As shown in \cite{cornwell2018existence}, an employed strategy can be used to compute $\text{Maslov}(w^*) $ for the doubly-diffusive FitzHugh-Nagumo system. Based on those, Definition \ref{def:another maslov index} is more practical than Definition \ref{def:maslov index for traveling wave} if one wants to compute the Maslov index for this case. Moreover, from the same discussion in the proof of Proposition \ref{pro:change index}, under (H\ref{as:4}), Definitions \ref{def:another maslov index} and  \ref{def:maslov index for traveling wave} are equivalent. The strategy employed in \cite{cornwell2018existence} may be valid for computing $\igeo(w^*)$ for the doubly diffusive FitzHugh-Nagumo system, such as to consider a traveling front solution established in \cite[Theorem 2.2]{cornwell2018existence}. 
\end{rem}

For the following FitzHugh-Nagumo equations 	\begin{align}\label{eq:f.n.eq.}
\begin{cases}
u_t&=u_{xx}+\frac{1}{d}(f(u)-v),\\
v_t&=v_{xx}+u-\gamma v,
\end{cases}
\end{align}
where $d,\  \gamma>0$ and $f(u)=u(1-u)(u-a)$ with $0<a<\frac{1}{2}$. If $\gamma $ is large enough, then the equation $u=\gamma f(u)$ has three solutions $0=u_1<u_2<u_3$, and we remark that it is easy to check that $f'(0)<0$ and $f'(u_3)<0$. Now, we consider a traveling front solution $ w^*$ of the FitzHugh-Nagumo equation \eqref{eq:f.n.eq.}: its
existence has been established in \cite{chen2018front}, and we remark that such waves are obtained as local minimizers of an energy functional, and based on the variational characterization, authors \cite{chenindex} utilize the spectral flow to define and calculate a stability index. As a supplement, we provide a geometric insight into the stability of the traveling front solution for the FitzHugh-Nagumo equation \eqref{eq:f.n.eq.}. 

Making a comparison with \eqref{eq:r.e.q.}, then, \eqref{eq:f.n.eq.} can be expressed as
\begin{align}
w_t=w_{xx}+QD\nabla F(w),
\end{align}where $D=\begin{bmatrix}
\frac{1}{d}&0\\0&1
\end{bmatrix}$, $Q=\begin{bmatrix}
1&0\\0&-1
\end{bmatrix}$ and $\nabla F(w)=\begin{bmatrix}
f(u)-v\\\gamma v-u
\end{bmatrix}.$
The eigenvalue problem \[\ddot{\phi}+c\dot{\phi}+QD\nabla^2F(w^*)\phi=\lambda \phi \]or its equivalent eigenvalue problem
\begin{align}
L\phi=\lambda \phi
\end{align}
is studied to determine the stability of $w^*$, where $L=\frac{\dif ^2}{\dif \xi^2}+c\frac{\dif }{\dif \xi}+QB(\xi)$ and $B(\xi)=D^\frac{1}{2}\nabla^2 F(w^*)D^{\frac{1}{2}}$.

\begin{thm}\label{thm:f.n.eq} Letting $w^*$ be a traveling front solution of the FitzHugh-Nagumo equation \eqref{eq:f.n.eq.} satisfying the following asymptotic condition: \[\lim\limits_{\xi\to-\infty}w^*=(u_3,\frac{u_3}{\gamma}) ,\ \lim\limits_{\xi\to+\infty}w^*=(0,0) \] and 
 $d>\gamma ^{-2}$, then we have that  \begin{equation}
	N_+(L)=\igeo(u^*).
	\end{equation}
\end{thm}


\section{Preliminary: an index formula }

As is commonly done for the Sturm-Liouville operators \cite{chen2007maslova,hu2020morse}, in this section, by applying \cite[Theorem 1]{hu2017index}, we derive an index formula (Cf. Proposition \ref{cor:spectral flow}) for the Hamiltonian system \eqref{eq:hamiltonian system}, which plays a crucial role in obtaining the spectral information of $L$. 

The notation of {\em  spectral flow \/}was introduced by Atiyah, Patodi and Singer in their study of index theory on manifolds with a boundary \cite{atiyah1976spectral}. For the reader's convenience, we provide a brief description of the basic properties of spectral flow.
Suppose that $E$ is a real separable Hilbert space, and  denote by $\cfsa(E)$ the space of all closed self-adjoint and Fredholm operators equipped with the {\em gap topology\/}. Let $ A:[0,1]\rightarrow\cfsa\left(E\right)$
be a continuous curve.
The $\Sf{A_t;t\in[0,1}]$ counts the algebraic multiplicities of the spectral flow of $A_t$ across the line $t=-\epsilon$ with some small positive number $\epsilon$.

For each $A_t$, there is an orthogonal splitting \[ E=E_-(A_t)\oplus E_0(A_t)\oplus E_+(A_t) ,  \] where $E_0$ is the null space of $A_t$, $\langle A_t v,v \rangle \geq 0 $ if $v\in E_+$ and $\langle A_t v,v \rangle \leq 0 $ if $v\in E_-$. 
There is an efficient way to compute the spectral flow through what are called 
{\em crossing 
	forms\/}. Let $P_t$ be the orthogonal projector from $E$ to 
$E_0\big( A_t\big)$. When 
$E_0\big(  A_{t_0}\big)\neq \{0\}$, we call
the instant $t_0$ a crossing instant. In this case, we defined the {\em crossing 
	from \/} 
$\mathrm{Cr}[A_{t_0}]$ as 
\[
\mathrm{Cr}[A_{t_0}]\= P_{t_0} \dfrac{\partial}{\partial t}P_{t_0}: 
E_0\big(  A_{t_0}\big)\to   
E_0\big(  A_{t_0}\big).
\]
We call the crossing instant $t_0$ {\em regular\/} if the crossing form 
$\mathrm{Cr}( A_{t_0})$ is 
nondegenerate. In this case we define the {\em signature\/} simply as 
\[
\sgn\big(\mathrm{Cr}(A_{t_0})\big)\= \dim E_+\big(\mathrm{Cr}(A_{t_0})\big)- 
\dim E_-\big(\mathrm{Cr}(A_{t_0})\big).
\]
We assume that all of the crossings are regular. Then, the crossing instants are isolated (and, hence, those on a compact interval are finite in number), and the spectral flow is given by the following formula 
\begin{equation}\label{eq:formula-sf-crossings}
\Sf{ A_t; t \in [0,1]}= \sum_{t_0 \in \mathcal S_*} 
\sgn\big(\mathrm{Cr}(A_{t_0})\big)- 
\dim E_-\big(\mathrm{Cr}(A_{0})\big)
+ \dim E_+\big(\mathrm{Cr}(A_{1})\big)
\end{equation}
where $\mathcal S_*\= \mathcal S\cap (a,b)$, and $\mathcal S$ denotes the set of 
all crossings.

 We list some basic properties of spectral flow:\begin{prop}\label{prop:properties of sf}
	\begin{enumerate}
		\item[(1)] If $t_0\in[0,1]$, then \[\Sf{A_t;t\in[0,1]}=\Sf{A_t;t\in[0,t_0]}+\Sf{A_t;t\in[t_0,1]}.\]\item[(2)] If $\mathcal S_*\= \mathcal S\cap (0,1)$ and $\mathcal S$ denotes the set of 
		all crossings, then\[  \Sf{A_t;t\in[0,1]}\leq \sum\limits_{t\in \mathcal S_* }\dim E_0( A_t).    \] 
	\end{enumerate}
\end{prop}


Now, we return to the problem of stability analysis.
If there is not a Hamiltonian structure for $L$ by the presence of the $\frac{\dif}{\dif \xi}$ term, we can circumvent this by considering instead the operator
\begin{equation}\label{eq:similar operator}
\mathbb{L}\=e^{\frac{c\xi}{2}}Le^{-\frac{c\xi}{2}}=\frac{\dif^2 }{\dif\xi^2}-\frac{c^2}{4}\Id+Q  B,
\end{equation}as is done in \cite{howard2016maslov,cornwell2019openinga}.
We now proceed to the study of the following eigenvalue problem \begin{equation}\label{eq:similar eigen.pro.}
\mathbb{L}\varphi =\lambda \varphi ,
\end{equation}
Setting $\psi= \dot \varphi-\frac{1}{2}c \varphi$ and $z=\begin{bmatrix}
\psi\\ \varphi
\end{bmatrix}$ converts Equation \eqref{eq:similar eigen.pro.} to the following Hamiltonian system \begin{equation}\label{eq:hamiltonian system}
\dot z=JH_\lambda(\xi)z,
\end{equation}where $H_\lambda=-JA_\lambda-\frac{1}{2}cJ=\begin{bmatrix}
Q &\frac{1}{2}cQ\\\frac{1}{2}cQ& B-\lambda Q 
\end{bmatrix}.$ We define the self-adjoint operators\begin{align}
 {F}_\lambda\=-J\frac{\dif}{\dif \xi}-H_\lambda:\dom {F}_\lambda\subset L^2(\R,\R^{2n})\subset L^2(\R,\R^{2n}),
\end{align}where $\dom  {F}_\lambda=W^{1,2}(\R,\R^{2n})$ and $\lambda\in[0,C]$. Here, $C$ is given in \eqref{eq:constant C}. 
The following Proposition \ref{lem:solution equivalent} shows that the Hamiltonian system \eqref{eq:hamiltonian system} has close contact with the system \eqref{eq:ODE}.

\begin{prop}\label{lem:solution equivalent}
	Letting $ \lambda\in \overline{\Com^+}$ and $z\in H^2(\R,\Com^{2n})$, then $z\in\ker\left(\frac{\dif}{\dif \xi}-A_\lambda \right)$ if and only if $e^{\frac{1}{2}c\xi}z\in\ker\left(-J\frac{\dif}{\dif \xi}-H_\lambda\right) $.
\end{prop}
It is well-known \cite[Lemma 3.1.10]{kapitula2013spectral} that the essential spectrum is given by \[\sigma_{ess}(L)=\{\lambda\in \Com|A_\lambda(+\infty)\mathrm{\ or\  }A_\lambda(-\infty) \mathrm{\ has\  a\ pure\ imaginary\ eigenvalue}\}.\]
  Following the same method as that in paper \cite[Section 2]{cornwell2018existence}, we here give the proof of the following Lemma for convenience of the reader.

\begin{lem}\label{lem:eigenvalue order of A}
	Under the condition (H\ref{as:1}), we have that \begin{itemize}
		\item[(1)] $\sigma_{ess}(L)\subset \Com^-$,
		\item[(2)] if $\lambda \in \overline{\Com^+}$, then $A_\lambda(\pm\infty)$ are hyperbolic.
	\end{itemize}
\end{lem}
\begin{proof}
	A simple calculation shows the eigenvalues of $A_\lambda(+\infty)$ and $A_\lambda(-\infty)$, such that \begin{equation}
	\mu(\lambda)=\frac{1}{2}\{-c\pm\sqrt{c^2+4(\lambda-\alpha)}\}\textrm{ and }\nu(\lambda)=\frac{1}{2}\{-c\pm\sqrt{c^2+4(\lambda-\beta),}\}\textrm{ respectively},
	\end{equation}where $\alpha$ is the eigenvalue of $Q B_+$, and $\beta$ is the eigenvalue of $Q B_-$. To prove that
	(1) holds, we only need to show that $A_\lambda(+\infty)$ and $A_\lambda(-\infty)$ have no purely imaginary
	eigenvalues if $\lambda\in \overline{\Com^+}$, which is equivalent to showing that $\re \sqrt{c^2+4(\lambda-\alpha)}\neq -c $ and $\re \sqrt{c^2+4(\lambda-\beta)}\neq -c $. From the formula $\re\sqrt{a+ib}=\{\frac{1}{2}\sqrt{a^{2}+b^{2}}+a\}^{\frac{1}{2}},$ we have that \begin{equation}\label{eq:real part}
	\re \sqrt{c^2+4(\lambda-\alpha)}=\{\frac{1}{2}[{\left(c^2+4\re\left(\lambda-\alpha\right)\right)^{2}+\left(\imaginary (\lambda-\alpha) \right)^{2}}]^{\frac{1}{2}}+c^2+4\re\left(\lambda-\alpha\right)\}^{\frac{1}{2}},
	\end{equation}we note that $\re\left(\lambda-\alpha\right)>0$, so Equation \eqref{eq:real part} implies that \begin{equation}
	\re \sqrt{c^2+4(\lambda-\alpha)}>c.
	\end{equation}
	This calculation actually proves that $A_\lambda(+\infty)$ has exactly $n$ eigenvalues of the positive real part and $n$ eigenvalues of the negative real part for $\re \lambda\geq 0$. Additionally, $A_\lambda(-\infty)$. We label the eigenvalues of $A_\lambda(\pm\infty)$ in order of increasing real part and observe that\begin{equation}
	\re \mu_1(\lambda)\leq \cdots \leq \re \mu_n(\lambda)<-c<0< \re \mu_{n+1}(\lambda)\leq \cdots \leq \re\mu_{2n}(\lambda),
	\end{equation}	
	\begin{equation}
	\re \nu_1(\lambda)\leq \cdots \leq \re \nu_n(\lambda)<-c<0< \re \nu_{n+1}(\lambda)\leq \cdots \leq \re\nu_{2n}(\lambda).
	\end{equation}From those, we complete the proof.
\end{proof}

\begin{rem}\label{rm:JB hyperbolic}
	We note that $JH_\lambda(\pm\infty)=A_\lambda(\pm\infty)+\frac{1}{2}c$: from  Lemma \ref{lem:eigenvalue order of A}, if  (H\ref{as:1}) holds, then $JH_\lambda(\pm\infty)$ are both hyperbolic for all $\lambda \in \overline{\Com^+}$, and this, together with \cite[Theorem 1]{hu2020morse}, tells us that for each $\lambda\in[0,C]$, the operator $ {F}_\lambda$ is a self-adjoint Fredholm operator, and in particular, it is possible to associate it to path \begin{align}\label{eq:F}
		\lambda\to   F_\lambda
	\end{align}the topological invariant: spectral flow. 
\end{rem}

\begin{proof}[The proof of Proposition \ref{lem:solution equivalent}] If $z\in \ker\left(\frac{\dif}{\dif \xi}-A_\lambda \right)$. By a simple calculation,\begin{align}\label{eq:solution = 1}
	\frac{\dif}{\dif \xi}\left( e^{\frac{1}{2}c\xi}z\right)=\frac{1}{2}ce^{\frac{1}{2}c\xi}z+e^{\frac{1}{2}c\xi}\dot z=J\left(-\frac{1}{2}cJ-JA_\lambda(\xi )\right)e^{\frac{1}{2}c\xi}z=JH_\lambda(\xi)e^{\frac{1}{2}c\xi}z.
	\end{align}    From Lemma \ref{lem:eigenvalue order of A}, $z $ must decay as fast as $e^{\nu_n(\lambda)\xi}$ as $\xi\to+\infty$, and we note that if $\nu_n(\lambda)<-c$, then $e^{\frac{1}{2}c\xi}z$ and $e^{\frac{1}{2}c\xi}\dot z$ both exponentially decay to $0$ as $\xi\to\pm\infty$, and hence, $e^{\frac{1}{2}c\xi}z\in L^2(\R,\Com^{2n})$ and $ e^{\frac{1}{2}c\xi}\dot z\in L^2(\R,\Com^{2n})$. Then, $e^{\frac{1}{2}c\xi}z\in H^2(\R,\Com^{2n})$: this together with Equation \eqref{eq:solution = 1} implies that 
	$e^{\frac{1}{2}c\xi}z\in\ker\left(-J\frac{\dif}{\dif \xi}-H_\lambda\right)$.
	
	Conversely, if $e^{\frac{1}{2}c\xi}z\in\ker\left(-J\frac{\dif}{\dif \xi}-H_\lambda\right)$, then by a simple calculation, \begin{align}\label{eq:solution = 2}
	\dot{z}=\frac{1}{2}cz+\dot z-\frac{1}{2}cz=e^{-\frac{1}{2}c\xi}\left(\frac{1}{2}ce^{\frac{1}{2}c\xi}z+e^{\frac{1}{2}c\xi}\dot z    \right)-\frac{1}{2}cz=e^{-\frac{1}{2}c\xi}\frac{\dif}{\dif \xi}\left(  e^{\frac{1}{2}c\xi} z  \right)-\frac{1}{2}cz=A_\lambda z.
	\end{align} 
		We note that when $JH_\lambda(\pm\infty)=A_\lambda(\pm\infty)+\frac{1}{2}c$, for each $\hat \mu\in\sigma(JH_\lambda(+\infty))\cap \Com ^+$, invoking Lemma \ref{lem:eigenvalue equivalent}, it is easy to check that $\Re \hat{\mu}>\frac{1}{2}c$: 
	then, $e^{\frac{1}{2}c\xi}z$ must decay at least as fast as $e^{\frac{1}{2}c\xi}$ as $\xi\to-\infty$. We thus determine that $z$ and $\dot z$ both exponentially decay as $\xi\to\pm\infty$, and then, $z,\ \dot z\in L^2(\R,\Com^{2n})$, and this together with Equation \eqref{eq:solution = 2}, implies that $z\in \ker\left(\frac{\dif}{\dif \xi}-A_\lambda \right)$.
	
\end{proof}
As a direct consequence, the following result holds.

\begin{cor}\label{lem:eigenvalue equivalent}Letting $\lambda\in \overline{\Com^+}$ and $\phi \in H^2(\R,\Com^n)$, then $\phi \in \ker(L-\lambda\Id)$ if and only if  $e^{\frac{c\xi}{2}}\phi \in \ker(\mathbb{L}-\lambda\Id)$.
\end{cor}

\begin{rem}\label{rm:eigenvalue equivalent}
	 Under the condition (H\ref{as:1}), if $y$ satisfies \eqref{eq:ODE} with $\lim\limits_{\xi\to+\infty}y=0$, then a similar discussion in the proof of Proposition \ref{lem:solution equivalent} guarantees that $e^{\frac{1}{2}\xi}y$ satisfies \eqref{eq:hamiltonian system} with $\lim\limits_{\xi\rightarrow+\infty}e^{\frac{1}{2}\xi}y=0$. Then, $\{e^{\frac{1}{2}c\tau}y(\tau)|y\text{ solve }(1.5) \text{ and } y\to0\text{ as }\tau\to+\infty\}$ of the Hamiltonian system \eqref{eq:hamiltonian system}, so we can say that $E^s_\lambda(\tau) $ is also the stable space of the Hamiltonian system \eqref{eq:hamiltonian system}. By the same reasoning, $E^u_\lambda(\tau)$ is also the unstable space of the Hamiltonian system \eqref{eq:hamiltonian system}. Moreover, letting $E_\lambda^s(+\infty)\=\left\{v\in \R^{2n}\left | \lim\limits_{\xi\to +\infty}\exp\left(\xi A_\lambda(+\infty)\right)v=0\right. \right\}\text{ and }E_\lambda^u(-\infty)\=\left\{v\in \R^{2n}\left | \lim\limits_{\xi\to -\infty}\exp\left(\xi A_\lambda(-\infty)\right)v=0\right. \right\}$, a similar discussion in Proposition \ref{lem:solution equivalent} shows that the following facts \\$E^s_\lambda(+\infty)=\left\{v\in \R^{2n}\left | \lim\limits_{\xi\to +\infty}\exp\left(\xi JH_\lambda(+\infty)\right)v=0\right. \right\}$ and $ E^u_\lambda(-\infty)=\left\{v\in \R^{2n}\left | \lim\limits_{\xi\to -\infty}\exp\left(\xi JH_\lambda(-\infty)\right)v=0\right. \right\} $ hold.
\end{rem} Under the condition (H\ref{as:1}), we determine that
   \begin{equation}
E^s_\lambda(+\infty)=\lim\limits_{\tau\to+\infty}E^s_\lambda(\tau) \text{ and } E^u_\lambda(-\infty)=\lim\limits_{\tau\to-\infty}E^u_\lambda(\tau),
\end{equation} 
where the convergence is meant in the gap (norm) topology of the Lagrangian Grassmannian (Cf. \cite{abbondandolo2003ordinary} for further details).

Given $\tau\geq 0$, let $B_1(\xi)=B(\tau+\xi) $ with $\xi\in\overline{\R^+}$ and $B_2(\xi)=B(\xi-\tau) $ with $\xi\in\overline{\R^-}$. The aim of the next part is to provide some sufficient condition for the coefficient of \eqref{eq:F} to obtain the nondegeneracy. 
\begin{lem}\label{lem:nondegenerate 2st ODE}
	Let \begin{equation}
	\mathbb{L}^+_{\lambda,M}\=\frac{\dif ^2}{\dif \xi^2}-\left\{\frac{c^2}{4}+\lambda\right\}\Id+Q B_1:W^{2,2}\left(\overline{\R^+},\R^n\right)\to L^2\left(\overline{\R^+},\R^n\right)
	\end{equation} and \begin{equation}
\mathbb{L}^-_{\lambda,M}\=\frac{\dif ^2}{\dif \xi^2}-\left\{\frac{c^2}{4}+\lambda\right\}\Id+Q B_2:W^{2,2}\left(\overline{\R^-},\R^n\right)\to L^2\left(\overline{\R^-},\R^n\right).
	\end{equation}
With $C$ given in \eqref{eq:constant C},	assuming that (H\ref{as:1}) holds, and $\lambda \geq C$, then the system\begin{align}
	\begin{cases}
\mathbb{L}^+_{\lambda,M}\varphi_1=0=\mathbb{L}^-_{\lambda,M}\varphi_2,\\
\varphi_1(0)=\varphi_2(0),\ \dot \varphi_1(0) =\dot \varphi_2(0) 
	\end{cases}
	\end{align}has only the zero solution.
\end{lem}
\begin{proof}
	Assuming that the system has a solution $(\varphi_1,\varphi_2)$, then we have that \begin{align}
	\langle\mathbb{L}^+_{\lambda,M}\varphi_1,\varphi_1\rangle_{L^2}+\langle\mathbb{L}^-_{\lambda,M}\varphi_2,\varphi_2\rangle_{L^2}=0
	\end{align}
	Integrating by part, we obtain
	\begin{align}
	&\langle \mathbb{L}^+_{\lambda,M} \varphi_1,\varphi_1  \rangle_{L^2}=-\norm{\dot \varphi_1}^2_{L^2}-\left\{\frac{c^2}{4}+\lambda\right\}\norm{\varphi_1}^2_{L^2}+\int_{0}^{+\infty}\langle Q B_1\varphi_1,\varphi_1\rangle\dif \xi-\langle \varphi_1(0),\dot \varphi_1(0)\rangle\\&\langle \mathbb{L}^-_{\lambda,M} \varphi_2,\varphi_2  \rangle_{L^2}=-\norm{\dot \varphi_2}^2_{L^2}-\left\{\frac{c^2}{4}+\lambda\right\}\norm{\varphi_2}^2_{L^2}+\int_{-\infty}^{0}\langle Q B_2\varphi_2,\varphi_2\rangle\dif \xi+\langle \varphi_2(0),\dot \varphi_2(0)\rangle
	\end{align}Let $I_{\lambda,1}\=-\norm{\dot \varphi_1}^2_{L^2}-\left\{\frac{c^2}{4}+\lambda\right\}\norm{\varphi_1}^2_{L^2}+\int_{0}^{+\infty}\langle Q B_1\varphi_1,\varphi_1\rangle\dif \xi$ and $I_{\lambda,2}\=-\norm{\dot \varphi_2}^2_{L^2}-\left\{\frac{c^2}{4}+\lambda\right\}\norm{\varphi_2}^2_{L^2}+\int_{-\infty}^{0}\langle Q B_2\varphi_2,\varphi_2\rangle\dif \xi$. It is easy to see that $I_i\leq -\norm{\dot \varphi_i}^2_{L^2}-\left\{\frac{c^2}{4}+\lambda-C\right\}\norm{\varphi_i}^2_{L^2} $. Then, by using the second condition in the above boundary value problem, we obtain
	\begin{equation}
	0=I_1+I_2\leq \sum\limits_{i=1,2} \left( -\norm{\dot \varphi_i}^2_{L^2}-\left\{\frac{c^2}{4}+\lambda-C\right\}\norm{\varphi_i}^2_{L^2}\right).
	\end{equation}If $\lambda \geq C$, then we infer that \[I_1+I_2=0\] if and only if $\varphi_i=\dot \varphi_i=0$ for $i=1,2$.  This concludes the proof.
	
\end{proof}

Let us now consider the associated first order differential operators $F_{\lambda,M}^+$ and $ F_{\lambda,M}^- $ of $\mathbb L^+_{\lambda,M}$ and $\mathbb L^-_{\lambda,M}$. A similar result holds.
\begin{lem}\label{lem:nondegenerate 1st ODE}
With $C$ given in \eqref{eq:constant C}, assuming that (H\ref{as:1}) holds, and $\lambda\geq C$, then the system\begin{align}
	\begin{cases}
	F_{\lambda,M}^+z_1=F_{\lambda,M}^-z_2=0\\
	z_1(0)=z_2(0)
	\end{cases}
	\end{align}has only the zero solution. 
\end{lem}
\begin{lem}\label{lem:non intersection lambda >>1}
With $C$ given in \eqref{eq:constant C},	assuming that (H\ref{as:1}) holds, we have that \begin{equation}
	E^s_\lambda(\tau)\cap E^u_\lambda (-\tau)=\{0\} \text{ for all }(\lambda,\tau)\in[C,+\infty)\times \overline{\R^+}.
	\end{equation}
\end{lem}
\begin{proof}
	Let $ B_1(xi)=  B(\xi+\tau)$ with $\xi\in\overline{\R^+}$ and $  B_2(\xi)=  B(\xi-\tau)$ with $\xi\in\overline{\R^-}$.  Then, the stable
	subspace of the equation $F_{\lambda,M}^+=0$ at $0$ is $E_\lambda^s$, the unstable
	subspace of the equation $F_{\lambda,M}^-=0$ at $0$ is $E_\lambda^u$, and there exists a linear bijection from the set of solutions
	of the system\begin{equation}
	\begin{cases}
		F_{\lambda,M}^+z_1=F_{\lambda,M}^-z_2=0\\
	z_1(0)=z_2(0)
	\end{cases}
	\end{equation}with the subspace $	E^s_\lambda(\tau)\cap E^u_\lambda (-\tau)$. By invoking once again Lemma \ref{lem:nondegenerate 1st ODE} and Lemma \ref{lem:nondegenerate 2st ODE}, we
	conclude that the initial value problem only admits the trivial solution for every $\lambda\geq C$. This concludes the proof.
	\end{proof}

By setting $ {B}_1(\xi)= {B}(\xi)$ for every $\xi\geq 0$ and $ {B}_2(\xi)= {B}(\xi)$ for every $\xi\leq 0 $, the following result holds.
\begin{lem}\label{lem:nondegenerate condition}
	With $C$ given in \eqref{eq:constant C}, assuming that (H\ref{as:1}) holds, if $\lambda \geq C$, then $\ker\left(\mathbb{L}-\lambda\right)=\{0\}$ and $\ker F_\lambda=\{0\}$.
\end{lem} From Lemma \ref{lem:nondegenerate condition}, we determine that the following result holds.
\begin{cor}\label{lem:nondegenerate condition S}
	With $C$ given in \eqref{eq:constant C}, assuming that (H\ref{as:1}) holds, and if $\lambda \geq C$, then we have that  $\ker S_\lambda =\{0\}$, where $S_\lambda\=-Q\frac{\dif^2}{\dif \xi^2}+\frac{c^2}{4}Q- {B}+\lambda  Q$.
\end{cor}It is well-known that for each $\lambda\in[0,C] $, the operator $S_\lambda$ is closed and self-adjoint with dense domain in $L^2(\R,\R^n )$. As a byproduct of condition (H\ref{as:1}) and \cite[Theorem 1]{hu2020morse}, $S_\lambda$ is also a Fredholm operator.

Finally, from \cite[Theorem 1]{hu2017index}, we obtain that \begin{align}\label{eq:spectral formula}
\Sf{F_\lambda,\lambda\in[0,C]}=&\iCLM(E^s_C(\tau),E^u_C(-\tau);\tau\in\overline{\R^+})-\iCLM(E^s(\tau),E^u(-\tau);\tau\in\overline{\R^+})\\&-\iCLM(E^s_\lambda(+\infty),E^u_\lambda(-\infty);\lambda\in[0,C]).
\end{align}
As a direct consequence of Lemma \ref{lem:non intersection lambda >>1} and Equation \eqref{eq:spectral formula}, we obtain the following result.
\begin{prop}\label{cor:spectral flow}
	Under the previous notations and assuming that (H\ref{as:1}) holds, the following equation holds:\begin{align}
	-\Sf{F_\lambda,\lambda\in[0,C]}=\iCLM(E^s(\tau),E^u(-\tau);\tau\in\overline{\R^+})+\iCLM(E^s_\lambda(+\infty),E^u_\lambda(-\infty);\lambda\in[0,C]).
	\end{align}
\end{prop}


\section{The proof of the {main result}}
The goal of this section is to provide a detailed proof of Theorem \ref{thm:main result}, Theorem \ref{thm:central result} and Theorem \ref{thm:f.n.eq}. Before showing the proof, we start by analyzing the distribution of the eigenvalues of $\mathbb{L}$.

As in \cite{chen2014stability}, the same discussion can be used here to obtain the distribution of eigenvalues of $\mathbb{L}$, which serves a crucial role of counting the number of nonnegative eigenvalues of $L$ via spectral flow. Let $Q^+$ and $Q^-$ be the orthogonal projections from $E$ to $E_+(Q)$ and $E_-(Q)$, respectively. Define $\mathscr S_1=Q^+SQ^+$, $\mathscr S_2=Q^-SQ^-$ and $\mathscr S_3=Q^+SQ^-$: in other words, $S$ can be decomposed as \begin{equation}
\begin{bmatrix}
\mathscr S_1&\mathscr S_3\\\mathscr S_3^*&\mathscr S_2
\end{bmatrix},
\end{equation}where $\mathscr S_3^*=\bar{\mathscr S}_3^T$ and $\bar{\mathscr S}_3$ denotes the complex conjugate of $\mathscr S_3$. For a linear self-adjoint operator $A$ defined on a Hilbert space $E$, denoted by $A>0$, if $\langle Av,v \rangle>$ for all $v\in E\backslash\{0\}$, two linear operators $A$ and $\hat A$ are denoted by $A>\hat A$ if $A-\hat A>0$.

\begin{lem}\cite{chen2014stability}\label{lem:positive eigenvalue}
	Supposing that $\mathscr S_1>0$ and $\Id>\mathscr S_3^*\mathscr S_1^{-2}\mathscr S_3$, then $\sigma(\mathbb{L})\cap \overline{{C}^+}\subset \R$. The same assertion holds if $-\mathscr S_2>0$ and $\Id >\mathscr S_3^*(-\mathscr S_2)^{-2}\mathscr S_3$.
\end{lem}

\begin{prop}\label{pro:sf one way}Under the condition (H\ref{as:1}), we have that \\
	(1) if $-\mathscr S_2>0$ and $\Id >\mathscr S_3^*(-\mathscr S_2)^{-2}\mathscr S_3$, then $\Sf{S_\lambda;\lambda\in[0,C]}=N_+(L)$,\\
	(2) if $\mathscr S_1>0$ and $\Id>\mathscr S_3^*\mathscr S_1^{-2}\mathscr S_3$, then $\Sf{S_\lambda;\lambda\in[0,C]}=-\overline N_+(L)$.
\end{prop}
\begin{proof}We only prove (1), while the other is analogous. 
	From Remark \ref{rm:eigenvalue equivalent}, we know that $N_+(L)=N_+(\mathbb{L})$. Next, we prove that $\Sf{S_\lambda;\lambda\in[0,C]}=N_+(\mathbb{L})$. Suppose that along the spectral flow, there is a crossing at $S_\lambda $ for some $\lambda \in[0,C]$ and $\phi\in \ker S_\lambda$: that is, \begin{align}\label{eq:above equation 1}
	-Q\ddot{\phi}+\left(\frac{1}{4}c^2Q-B+\lambda Q\right)\phi=0.
	\end{align}Letting $\phi_+=Q^+\phi$ and $\phi_-=Q^-\phi $, we can rewrite Equation \eqref{eq:above equation 1} as \begin{align}\label{eq:above equation 2}
	\mathscr S_1\phi_++\mathscr S_3\phi_-=-\lambda\phi_+\end{align}
	\begin{align}\label{eq:above equation 3}
	\mathscr S^*_3\phi_++\mathscr S_2\phi_-=\lambda\phi_-.
	\end{align}

	Solving Equation \eqref{eq:above equation 3}, we obtain that $ \phi_-=(\lambda-\mathscr S_2)^{-1}\mathscr S_3^*\phi_+$, and this, together with Equation \eqref{eq:above equation 2}, obtains \begin{align}
	\frac{\dif }{\dif \lambda }\langle S_\lambda \phi,\phi \rangle =\langle Q\phi,\phi\rangle=\langle \phi_+,\phi_+\rangle-\langle \phi_-,\phi_-\rangle=\langle \phi_+,\phi_+\rangle-\langle \mathscr S_3^*(\lambda-\mathscr S_2)^{-2}\mathscr S_3\phi_+,\phi_+  \rangle 
	\end{align}We note that $\Id >\mathscr S_3^*(-\mathscr S_2)^{-2}\mathscr S_3>\mathscr S_3^*(\lambda-\mathscr S_2)^{-2}\mathscr S_3$ for all $\lambda \geq 0$. This indicates that the sign
	of the crossing form has to be positive whenever a crossing occurs at $\lambda \in[0,C ]$. In view of Equation \eqref{eq:formula-sf-crossings}, we conclude from Lemma that \begin{equation}
	\Sf{S_\lambda;\lambda\in[0,C]}=\sum\limits_{\lambda\in(0,C]}\dim \ker S_\lambda.
	\end{equation}For (2), a slightly modified argument shows that the sign of the crossing operator must be negative if a crossing occurs at $\lambda\in[0,C]$, and then \begin{equation}
	\Sf{S_\lambda;\lambda\in[0,C]}=-\sum\limits_{\lambda\in[0,C]}\dim \ker S_\lambda.
	\end{equation}This completes the proof.
\end{proof}

The aim of the next part is to prove some transversal properties about some invariant subspaces that are useful in our proof.

\begin{lem}\label{lem:transiversal result}
	Under the conditions (H\ref{as:1}) and (H\ref{as:3}), we have that \begin{equation}
	E^s_\lambda(+\infty)\pitchfork L_R \text{ and } E^u_\lambda(-\infty)\pitchfork L_R,
	\end{equation}.
\end{lem}
\begin{proof}
	We provide the proof of $E^s_\lambda(+\infty)\pitchfork L_R$ in completely similar fashion. Let $\begin{bmatrix}
	p\\q
	\end{bmatrix}\in E^s_\lambda(+\infty)\cap  L_R$, and noting that $E^s_\lambda(+\infty)$ is invariant under $JH_\lambda(+\infty)$, then $JH_\lambda(+\infty)\begin{bmatrix}
	p\\q
	\end{bmatrix}\in E^s_\lambda(+\infty)$. From (H\ref{as:3}), a direct computation yields that \begin{align}
	0&=\omega\left(JH_\lambda(+\infty)\begin{bmatrix}
	p\\q
	\end{bmatrix},\begin{bmatrix}
	p\\q
	\end{bmatrix}\right)=-\Inp{\begin{bmatrix}
		Q&\frac{1}{2}cQ\\\frac{1}{2}cQ& B_+-\lambda Q^{-1}
		\end{bmatrix}\begin{bmatrix}
		p\\q
		\end{bmatrix},\begin{bmatrix}
		p\\q
		\end{bmatrix}}\\&=-\Inp{\begin{bmatrix}
		p-\frac{1}{2}cq\\\frac{1}{2}cp+\lambda q+ B_+q
		\end{bmatrix},\begin{bmatrix}
		p\\q
		\end{bmatrix}}=-  \Inp{p,p}  -\lambda \Inp{q,q}  +\Inp{Q B_+ q,q}\leq 0,
	\end{align} we determine that $E^s_\lambda(+\infty)\pitchfork L_R$. This completes the proof.
\end{proof}
Now, following Lemma \ref{lem: maslov triple index (pair)} and Equation \eqref{eq:the triple index}, we have that \begin{align}\label{eq:compute maslov index}
\iCLM\left(E^s_\lambda(+\infty),E^u_\lambda(-\infty);\lambda\in[0,C]\right)=&\igeo\left(E^u_C(-\infty),E^s_C(+\infty);L_R\right)-\igeo\left(E^u (-\infty),E^s (+\infty);L_R\right)\\=&m^+(\mathfrak{Q}\left(E_C^u(-\infty),E_C^s(+\infty);L_R\right)) -m^+(\mathfrak{Q}\left(E^u(-\infty),E^s(+\infty);L_R\right)) .
\end{align}

We recall that a Lagrangian frame for a Lagrangian subspace $L$ is an injective linear map $T:\R^n\to \R^{2n}$ whose image is $L$ (Cf. page 828 of \cite{robbin1993maslov}). Such a frame has the form \[T=\begin{bmatrix}
X\\Y
\end{bmatrix},\] where $X$ and $Y$ are both $n\times n$ matrices and $X^TQY=Y^TQX$.  

We introduce the following notations $T_\lambda^+=\begin{bmatrix}
X_{\lambda,+}&Y_{\lambda,+}\\0&\Id_2\\\Id_1&0\\Y^T_{\lambda,+}&Z_{\lambda,+}
\end{bmatrix}$, $T_\lambda^-=\begin{bmatrix}
X_{\lambda,-}&Y_{\lambda,-}\\0&\Id_2\\\Id_1&0\\Y^T_{\lambda,-}&Z_{\lambda,-}
\end{bmatrix}$, $M_{\lambda,+}=\begin{bmatrix}
X_{\lambda,+}&Y^T_{\lambda,+}\\Y_{\lambda,+}&Z_{\lambda,+}
\end{bmatrix}$ and $M_{\lambda,-}=\begin{bmatrix}
X_{\lambda,-}&Y_{\lambda,-}^T\\Y_{\lambda,-}&Z_{\lambda,-}
\end{bmatrix}$, where $X_{\lambda,\pm}$ are $r\times r$ matrices, $Y_{\lambda,\pm}$ are $ (n-r)\times r$ matrices, $Z_{\lambda,\pm}$ are $(n-r)\times (n-r)$ matrices, $\Id_1$ is the $r\times r$ identity matrix, and $I_2$ is the $(n-r)\times (n-r)$ identity matrix. Here, $r=\dim V^+(Q)$. For each $\lambda\in\overline{\R^+}$, from Lemma \ref{lem:transiversal result}, we can use notations $T_{\lambda,+}$ and $T_{\lambda,-}$ for the Lagrangian frames of $E^s_\lambda(+\infty)$ and $E^u_\lambda(-\infty)$, respectively.

We now consider the following operator: \begin{equation}
F^+_\lambda\=-J\frac{\dif}{\dif \xi}-H_\lambda(+\infty)
\end{equation} and the associated second order operator $\mathbb L_\lambda^+$ and let $\varphi$ be a solution of $\mathbb L_\lambda^+\varphi=0$, where $\mathbb L_{\lambda,M}^+$ denotes the operator $\mathbb L_{\lambda}^+$ defined on the maximal domain $W^{2,2}(\overline{\R^+}, \R^n)$. Then, the map $\phi\mapsto  (\dot \phi^T(0)-\frac{1}{2}c \phi^T(0),\phi^T(0))^T$ provides a linear bijection from $\ker \mathbb L_{\lambda,M}^+$ to $E^s_\lambda(+\infty)=V^-(JH_\lambda(+\infty))$.

We note that for each $z\in E^s_\lambda(+\infty)$, there exists  $u=\begin{bmatrix}
p\\q
\end{bmatrix}\in V^+(Q)\oplus V^-(Q) \cong\R^n$, such that $z=T_{\lambda,+}(u)$.

Let $\varphi(\xi)\in \ker \mathbb L_\lambda^+$ with $(\dot \phi^T(0)-\frac{1}{2}c \phi^T(0),\phi^T(0))^T=T^+_\lambda u\in E^s_\lambda(+\infty)$, where  $u=\begin{bmatrix}
p\\q
\end{bmatrix}\in V^+(Q)\oplus V^-(Q) \cong\R^n$. A simple calculation shows that \begin{align}
0&=\langle \mathbb L_\lambda^+ \varphi(\xi),\varphi(\xi)\rangle_{L^2}=\langle\dot \varphi,\dot \varphi\rangle_{L^2}-c\langle \varphi,\dot \varphi\rangle_{L^2}+\frac{c^2}{4}\norm{\varphi}^2_{L^2}+\lambda \norm{\varphi}^2_{L^2}-\langle Q B_+\varphi,\varphi\rangle_{L^2}+\langle \dot \varphi(0)-\frac{c}{2}\varphi(0),\varphi(0)\rangle\\
&= \norm{\dot \varphi-\frac{c}{2}\varphi}^2_{L^2}-\int_{0}^{+\infty}\Inp{QB_+\varphi,\varphi} \dif\xi+\lambda \norm{\varphi}^2_{L^2}+ \Inp{M_{\lambda,+}\begin{bmatrix}
	p\\q
	\end{bmatrix}, \begin{bmatrix}
	p\\q
	\end{bmatrix}} \\&\geq \norm{\dot \varphi-\frac{c}{2}\varphi}^2_{L^2}+(\lambda-C)\norm{\varphi}^2_{L^2}+\Inp{M_{\lambda,+}u, u},
\end{align}this equation, together with (H\ref{as:1}), (H\ref{as:3}) and (H\ref{as:4}), shows the following Lemma:
\begin{lem}\label{lem:matrix<0}
	With $C$ given in \eqref{eq:constant C}, the following results hold: \\(1) if (H\ref{as:1}) and (H\ref{as:3}) hold, we have that $M_{\lambda,+}(s)$ is negative definite for all $\lambda \geq  C$. \\(2) if (H\ref{as:4}) holds, we have that $M_{\lambda,+}(s)$ is negative definite for all $\lambda \geq  0$.
\end{lem}
Similarly, we have that 
\begin{lem}\label{lem:matrix>0}
With $C$ given in \eqref{eq:constant C},	 the following results hold: \\(1) if (H\ref{as:1}) and (H\ref{as:3}) hold, we have that $M_{\lambda,-}$ is positive definite for all $\lambda \geq  C$. \\(2) if (H\ref{as:4}) holds, we have that $M_{\lambda,-}$ is positive definite for all $\lambda \geq  0$.
\end{lem}
Letting $z\in E^u_\lambda(-\infty)$, then there is $u=\begin{bmatrix}
p\\q
\end{bmatrix}\in V^+(Q)\oplus V^-(Q) \cong\R^n$, such that $z=T_{\lambda,+}u$, and $z$ can be rewritten as $z=T_{\lambda,+}u+(T_{\lambda,-}-T_{\lambda,+})u\in E^s_\lambda+L_R$. 
From a simple calculation, we have that \begin{align}\label{eq:form Q}
&\mathfrak{Q}\left(E^u_\lambda(-\infty),E^s_\lambda(+\infty);L_R\right)(u,u)\\=&\omega\left(\begin{bmatrix}X_{\lambda,+}p+Y^T_{\lambda,+}q\\q\\p\\Y_{\lambda,+}p+Z_{\lambda,+}q
\end{bmatrix},\begin{bmatrix}
\left(X_{\lambda,-}-X_{\lambda,+}\right)p+\left(Y^T_{\lambda,-}-Y^T_{\lambda,+}\right)q\\0\\0\\\left(Y_{\lambda,-}-Y_{\lambda,+}\right)p+\left(Z_{\lambda,-}-Z_{\lambda,+}\right)q
\end{bmatrix}\right)\\=&\Inp{(X_{\lambda,+}-X_{\lambda,-})p,p}+2\Inp{(Y^T_{\lambda,+}-Y^T_{\lambda,-})p,q}+\Inp{(Z_{\lambda,+}-Z_{\lambda,-})q,q}\\=&\Inp{\begin{bmatrix}X_{\lambda,+}-X_{\lambda,-}&Y^T_{\lambda,+}-Y^T_{\lambda,-}\\Y_{\lambda,+}-Y_{\lambda,-}&Z_{\lambda,+}-Z_{\lambda,-}
	\end{bmatrix}\begin{bmatrix}
	p\\q
	\end{bmatrix},\begin{bmatrix}
	p\\q
	\end{bmatrix}}=\Inp{(M_{\lambda,+}-M_{\lambda,-})u,u}.
\end{align}
From Lemma \ref{lem:matrix<0}, Lemma \ref{lem:matrix>0}, Equation \eqref{eq:form Q} and Equation \eqref{eq:the triple index}, the following result holds.
\begin{lem}\label{lem:boundary maslov index = 0}
	With $C$ given in \eqref{eq:constant C}, the following results hold: \\(1) if  (H\ref{as:1}) and (H\ref{as:3}) hold, then we have that \begin{equation}
	\igeo\left(E^u_\lambda(-\infty),E^s_\lambda(+\infty),L_R\right)=0\ \text{for }\lambda\geq C.
	\end{equation}(2) if (H\ref{as:4}) holds, then we have that \begin{equation}
	\igeo\left(E^u_\lambda(-\infty),E^s_\lambda(+\infty),L_R\right)=0\ \text{for }\lambda\geq 0.
	\end{equation}
\end{lem}

Before finishing the preparation of our proof of our main results, we recall the definition of {\em positive curve\/}. 
\begin{defn}\cite{hu2020decomposition}\label{def:positive curve} 
	Let $ A:[0,1]\rightarrow\cfsa\left(E\right)$
	be a continuous curve. The curve $  A$ is named a {\em  positive curve\/} if $\Set{t|\ker  A_\lambda\neq 0}$ is finite and 
	\[
	\Sf{ A_t;t\in[0,1]}=\sum_{0<t\leq 1}\dim\ker  A_t  .
	\]
\end{defn}

\begin{proof}[\textbf{The Proof of Proposition \ref{pro:change index}}]


	For some $a\in\R$, we construct the following homotopy Lagrangian path\begin{align}
		\left(E^s(\tau+sa),E^u(-\tau+sa)\right),(\tau,s)\in\R^+\times [0,1].
	\end{align}
	We point out that $\dim\left(E^s(sa)\cap E^u(sa)\right)$ is constant for all $s\in[0,1]$ and $E^s(+\infty)\pitchfork E^u(-\infty)$.

	By the stratum homotopy invariance property of the Maslov index, we have that \begin{align}\label{eq:maslov index change}
		&\iCLM\left(E^s(\tau),E^u(-\tau);\tau\in\R^+     \right)=\iCLM\left(E^s(\tau+ {a}),E^u(-\tau+ {a});\tau\in\R^+     \right)\\=&\iCLM\left(E^s(\tau+2 {a}),E^u(-\tau);\tau\in[- {a},+\infty)     \right).
	\end{align} 
	We know that $E^s(+\infty)\pitchfork E^u(-\infty)$ and $E^s(\tau)\to E^s(+\infty)$ as $\tau\to+\infty$ under the gap topology of the Lagrangian Grassmannian, so we can choose $\tau_0$, such that $E^s(\tau)\pitchfork E^u(-\infty)$ for all $\tau\geq \tau_0$, and the path $E^u(\tau):(-\infty,\tau_0]\to \Lagr(n)$ has only regular crossing with respect to $E^s(\tau_0)$. Letting $ {a}=\tau_0$, we construct the following homotopy Lagrangian path: 
	\[ 
	\left(E^s(\tau_0+s(\tau_0+\tau)),E^u(-\tau)\right),\ (\tau,s)\in[-\tau_0,+\infty) \times[0,1].
	\]By the stratum homotopy invariance, reversal property of Maslov index and Equation \eqref{eq:maslov index regular crossing}, we determine that \begin{align}
		&\iCLM\left(E^s(\tau+2\tau_0),E^u(-\tau);\tau\in[-\tau_0,+\infty)    \right)=\iCLM \left(E^s(\tau_0),E^u(-\tau);\tau\in[-\tau_0,+\infty) \right)\\&=-\iCLM\left( E^s(\tau_0),E^u(\tau);\tau\in(-\infty,\tau_0]  \right)=-\mathrm{Maslov }(w^*),
	\end{align} and this, together with Equation \eqref{eq:maslov index change}, completes the proof.

\end{proof}

\begin{proof}[\textbf{The Proof of Theorem \ref{thm:main result}}]

 We first prove that $ \Sf{ S_\lambda;\lambda\in[0,C]}=\Sf{F_\lambda;\lambda\in[0,C]}$.
	We start by introducing the continuous map 
	\[
	f:\cfsa\left(L^2\left(\R,\R^{n}\right)\right)\rightarrow\cfsa\left(L^2\left(\R,\R^{2n}\right)\right) \textrm{ defined by } f( S_\lambda)\= F_\lambda.
	\]
	Let $h(\lambda ,s)=f( S_\lambda+s\Id)$ for $(\lambda,s)\in[0,C]\times[0,\epsilon]$. Then, for every $\lambda\in[0,C]$, $h(\lambda,s)$ is a positive curve.
	Let $\lambda_0 \in[0,C]$ be a crossing instant for the path $\lambda \mapsto  S_\lambda$, meaning that $\ker   S_{ \lambda_0}\neq \{0\}$, and let us consider the positive path $s \mapsto S_{ \lambda_0 }+s\Id$. Thus, there exists $\delta>0$, such that $\ker\left(  S_{\lambda_0} +\delta\Id \right)=\{0\}$, which is equivalent to  $\ker h(\lambda_0,\delta)=\{0\}$. Since $S_{\lambda_0}+\delta\Id$ is a Fredholm operator, then there exists $\delta_1>0$, such that $\ker(S_\lambda+\delta\Id)=\{0\}$ for every $\lambda\in[\lambda_0-\delta_1,\lambda_0+\delta_1]$. By this argument, we determine that $ \ker h(\lambda,\delta)=\{0\}$ for every  $\lambda \in[\lambda_0-\delta_1,\lambda_0+\delta_1]$. Then, we determine that 
	\begin{equation}
	\begin{cases}
	\Sf{S_\lambda+\delta\Id,\lambda\in[\lambda_0-\delta_1,\lambda_0+\delta_1]}=0,\\
	\Sf{ h(\lambda,\delta), \lambda\in[\lambda_0-\delta_1,\lambda_0+\delta_1]}=0.
	\end{cases}
	\end{equation}   
	By the homotopy invariance of the spectral flow, we infer that 
	\begin{equation}\label{eq:spectral flow A}
	\Sf{S_\lambda,\lambda\in[\lambda_0-\delta_1,\lambda_0+\delta_1]}=\Sf{S_{\lambda_0-\delta_1}+s\Id,s\in[0,\delta]}-\Sf{S_{\lambda_0+\delta_1}+s\Id,s\in[0,\delta]}
	\end{equation}
	and 
	\begin{equation}\label{eq:spectral flow F}
	\Sf{h(\lambda,0),\lambda\in[\lambda_0-\delta_1,\lambda_0+\delta_1]}=\Sf{h(\lambda_0-\delta_1,s),s\in[0,\delta]}-\Sf{h(\lambda_0-\delta_1,s),s\in[0,\delta]}
	\end{equation} 
	We observe that $s\mapsto S_{\lambda_0\pm\delta_1}+s\Id $ and $s\mapsto h(\lambda_0\pm\delta_1,s)$ are both positive curves. It follows that 
	\begin{align}\label{eq:sf A= sf F}
	\Sf{ S_{\lambda_0\pm\delta_1}+s\Id,s\in[0,\delta]}&=\sum_{0<s\leq\delta}\dim\ker\left( S_{\lambda_0\pm\delta_1}+s\Id \right)=\sum_{0<s\leq\delta}\dim\ker h(\lambda_0\pm\delta_1,s) \\
	&=\Sf{h(\lambda_0\pm\delta_1,s),s\in[0,\delta]} 
	\end{align}
	  from Equations~\eqref{eq:spectral flow A}, \eqref{eq:spectral flow F} and \eqref{eq:sf A= sf F}, we have that \begin{align}\label{eq:sf = local}
		&\Sf{S_\lambda,\lambda\in[\lambda_0-\delta_1,\lambda_0+\delta_1]}=\Sf{h(\lambda,0),\lambda\in[\lambda_0-\delta_1,\lambda_0+\delta_1]}\\=&\Sf{{F}_\lambda,\lambda\in[\lambda_0-\delta_1,\lambda_0+\delta_1]}.
		\end{align} the crossing instants are isolated, and those on a compact interval are finite in number. From Equation~\eqref{eq:sf = local} and the path additivity of spectral flow, we determine that  $\Sf{S_\lambda;\lambda\in[0,C]}=\Sf{F_\lambda;\lambda\in[0,C]}$.

		From Proposition \ref{cor:spectral flow} and Equation \eqref{eq:compute maslov index}, we determine that 
	 \begin{align}\label{eq:spectral flow}
-\Sf{S_\lambda;\lambda\in[0,C]}=\iCLM\left(E^s(\tau),E^u(-\tau);\tau\in \overline{\R^+}\right)-\igeo\left(E^u(-\infty),E^s(+\infty);L_R\right).
		\end{align} 
	Since (H\ref{as:1}) and (H\ref{as:3}) hold, from Lemma \ref{lem:boundary maslov index = 0}, Equation \eqref{eq:spectral flow} and Proposition \ref{prop:properties of sf}, then \begin{equation}
	|\igeo(w^*)+\igeo\left(E^u(-\infty),E^s(+\infty);L_R\right)|\leq \overline N_+(L),
	\end{equation}
\end{proof}
\begin{proof}[\textbf{The proof of Theorem \ref{thm:central result}}]
	From Remark \ref{rem:condition}, (H\ref{as:4}) implies that (H\ref{as:1}) holds, and then, from Theorem \ref{thm:main result} and Lemma \ref{lem:boundary maslov index = 0}, we have that \begin{equation}
		|\igeo(w^*)|\leq\overline  N_+(L).
	\end{equation} 
\end{proof}

\begin{proof}[\textbf{The proof of Theorem \ref{thm:f.n.eq}}]
	We note that $S=\begin{bmatrix}
	-\frac{\dif^2}{\dif t^2}+\frac{c^2}{4}-\frac{f'(u^*)}{d}&\frac{1}{\sqrt{d}}\\\frac{1}{\sqrt{d}}&\frac{\dif^2}{\dif t^2}-\frac{c^2}{4}-\gamma 
	\end{bmatrix}$, and then, $\mathscr S_3=\frac{1}{\sqrt{d}}$, $\mathscr S_2= \frac{\dif^2}{\dif t^2}-\frac{c^2}{4}-\gamma$. Since $\gamma>0$, it is easy to see that $\mathscr S_2>0$. Moreover, if $d>\gamma^{-2}$, then we have that $\Id>\frac{1}{d}\left(-\frac{\dif^2}{\dif t^2}+\frac{c^2}{4}+\gamma \right)^{-2}=\mathscr S^*_3(-\mathscr S_2)^{-2}\mathscr S_3$, so Proposition \ref{pro:sf one way} (1) holds, and then, we have that $\Sf{S_\lambda;\lambda\in[0,C]}=N_+(L) $. Moreover, by a simple calculation 
	and the facts $f'(0)<0$ and $f'(u_3)<0$, it is easy to check that the condition   (H\ref{as:4}) holds. Then, from Equation \eqref{eq:spectral flow} and Lemma \ref{lem:boundary maslov index = 0}, we complete the proof.
\end{proof}


%
%
%
%




\section{The triple and H\"ormander index}\label{subsec:Maslov}

Recently, Zhu et al., in the interesting  paper \cite{zhou2018hormander}, deeply investigated the H\"ormander index, particularly its relation with respect to the so-called {\em triple index\/} in a slightly generalized (in fact, isotropic) setting. Given three isotropic subspaces  $\alpha,\beta$ and $\delta$ in  $(\R^{2n},\omega)$, we define the quadratic form $\mathfrak Q$ as follows: 
\begin{equation}\label{eq:the form Q}
\mathfrak{Q}\=\mathfrak{Q}(\alpha,\beta;\delta): \alpha\cap(\beta+\delta)\to \R \quad \textrm{ given by } \quad 
\mathfrak{Q}(x_1,x_2)=\omega(y_1,z_2), 
\end{equation}
where for $j=1,2$, $x_{j}=y_{j}+z_{j}\in\alpha\cap(\beta+\delta)$ and $y_{j}\in\beta$, $z_{j}\in\delta$. By invoking \cite[Lemma 3.3]{zhou2018hormander}, in the particular case in which $\alpha,\beta,\delta$ are Lagrangian subspaces, we obtain
\begin{equation}\label{Qker}
\ker\mathfrak{Q}(\alpha,\beta;\delta)=\alpha\cap\beta+\alpha\cap\delta.
\end{equation}
By \cite[Lemma 3.13]{zhou2018hormander}, we are in position to define the {\em triple index\/} in terms of the quadratic form  $\mathfrak Q$ defined above.
\begin{defn}\label{eq:trip1} Let $\alpha,\beta$ and $\kappa$ be	three Lagrangian subspaces of symplectic vector space $(\R^{2n},\omega)$. Then, the {\em triple index of the triple $(\alpha,\beta,\kappa)$\/} is defined by 
	\begin{equation}\label{eq:the triple index}
		\iota(\alpha,\beta,\kappa)=m^{+}\big(\mathfrak{Q}(\alpha,\beta;\kappa)\big)+\dim\big(\alpha\cap\kappa)-\dim(\alpha\cap\beta\cap\kappa)\big).
	\end{equation}
	where $m^+$ is the Morse positive index of a quadratic form Q. 
\end{defn}

Another closely related symplectic invariant is the so-called {\em  H\"{o}rmander index\/}, which is particularly important for measuring the difference in the (relative) Maslov index computed with respect to two different Lagrangian subspaces (we refer the interested reader to the celebrated and beautiful paper \cite{robbin1993maslov} and the references therein). 

Let $V_0,V_1,L _0,L _1$ be four Lagrangian subspaces and 
$L\in \mathscr C^0\big([0,1],\Lagr(n)\big)$ be such that $L(0)=L_0$ and $L(1)=L_1$. 
\begin{defn}\label{def:hormander index} 
	Letting $L,\ V\in\mathscr C^0([0,1],\mathrm{Lag}(n))$ be such that $L(0)=L_0,\ L(1)=L_1,\ V(0)=V_0$ and $V(1)=V_1$, the H\"ormander index is the integer defined by \begin{align}
	s(L_0,L_1;V_0,V_1)&=\iCLM(V_1,L(t);t\in[0,1])-\iCLM(V_0,L(t));t\in[0,1]\\&=\iCLM(V(t),L_1;t\in[0,1])-\iCLM(V(t),L_0;t\in[0,1])
	\end{align}
\end{defn}
\begin{rem}
	As a direct consequence of the fixed endpoints homotopy invariance of the $\iCLM$-index, it is actually possible to prove that Definition~\ref{def:hormander index} is well-posed, meaning that it is independent of the path $L$ joining the two Lagrangian subspaces $L_0,L_1$. (Cf. \cite{robbin1993maslov} for further details). 
\end{rem}

Let us now be given four Lagrangian subspaces, namely $\lambda_1,\lambda_2,\kappa_1,\kappa_2$ of symplectic vector space $(\R^{2n},\omega)$. By \cite[Theorem 1.1]{zhou2018hormander}, the H\"{o}rmander index $s(\lambda_1,\lambda_2;\kappa_1,\kappa_2)$ can be expressed in terms of the triple index as follows 
\begin{equation}\label{eq:the Hormander index computed by triple index}
s(\lambda_1,\lambda_2;\kappa_1,\kappa_2)=\iota(\lambda_1,\lambda_2,\kappa_2)-\iota(\lambda_1,\lambda_2,\kappa_1)=\iota(\lambda_1,\kappa_1,\kappa_2)-\iota(\lambda_2,\kappa_1,\kappa_2).
\end{equation}

\begin{lem}\cite{hu2020morse}	\label{lem: maslov triple index (pair)} 
	Let $L_1(t)$ and
	$L_2(t)$ be two paths in $\Lagr(n)$  with
	$t\in[0,1]$, and assume that $L_1(t)$ and $L_2(t)$
	are both transversal to the (fixed) Lagrangian subspace $L$. We then obtain  
	\[
	\iCLM\big(L_1(t),L_2(t);t\in[0,1]\big) =\iota\big(L_2(1),L_1(1);L)-\iota(L_2(0),L_1(0);L).
	\]
\end{lem}

\section*{Acknowledgments} The author is grateful to Professors Xijun Hu and Li Wu for their helpful discussion during the preparation of this article.
Finally, the author thanks the anonymous referee for reading the paper carefully and providing thoughtful comments.

%
%
%
%

\printbibliography

\end{document}